\documentclass[12pt,ngerman]{article}

\usepackage{amsmath,amsthm,amssymb,amscd}
\usepackage[ansinew]{inputenc}
\usepackage{graphicx}
\usepackage{color}

\begin{document}
\title{A note on reassignment for wavelets
}

\author{Hans Martin Reimann, Berne}
\date{28th of September, 2015}
\maketitle

\vspace{4cm}
Abstract\\
The reassignment method for the wavelet transform is investigated. Particularly good results are obtained if the wavelet is an extremal for the uncertainty relation of the affine group.
\vspace{4cm}

\section{Introduction}

Reassignment was introduced 1976 by Kodera, De Villedary and Gendrin [KDG] as a numerical method for improving the readabilty of time-frequency representations for nonstationary signals. Its main applicability is in the analysis of spectograms of the short time Fourier transform. In their paper of 2006, Gardner and Magnasco make the point that the auditory system uses sparse time-frequency representations and they hold that hearing uses reassignment or a homologous method for time-frequency processing. It is not known how exactly the brain processes temporal and spectoral information, yet reassignment displays some of the essential features which seem to be inherent in the way that the brain is processing signals: it uses simple operations (on phase and - as will become clear - on amplitude), it makes implicit use of time and frequency (the latter will here be

replaced by scale) and it uses a multiplicity of bandwidths.\\
In this note the reassignment method for the short time Fourier transform is reviewed. It is then shown that the reassignment method for the wavelet transform is apt to lead to good results, provided the wavelet is appropriately chosen. In the context of the windowed Fourier transform, the Gaussian function plays a particularly prominent role. Gaussians are extremal functions for the Heisenberg uncertainty role and the windowed Fourier transform with a Gaussian is related to the representations of the Heisenberg group. The wavelet transformation goes with the affine group. It is an intertwininig transformation between two of its representations. The specific wavelets that show optimal behaviour are the extremals for the uncertainty relation for the affine group. In [R] it is shown that they are related to the auditory process. This is then the main motivation to study the reassignment method in the context of the wavelet transform. It is argued in this note that the wavelet transform with wavelets that are extremal for the uncertainty relation for the affine group behave particularly well under reassignment.

\vspace{4cm}
\section{Short time Fourier transform}
The short time Fourier transform with window $h$ is defined as
\begin{eqnarray*}
Sf(t,\omega)&=&
  \int f(s)\bar{h}(s-t)e^{-i\omega (s-\frac{t}{2}) }ds\\
&=& e^{i\frac{\omega t}{2}}  \int f(s)\bar{h}(s-t)e^{-i\omega s }ds
\end{eqnarray*}

The so-called Weyl factor $  e^{i\frac{\omega t}{2}}$ connects this transform to the Schroedinger representation of the Heisenberg group. The Heisenberg group is the space $\mathbb{R}^3 = \{(x=(s,\xi, z)\}$ equipped with the multiplication law
\begin{equation}
x_1. x_2 = (s_1 +s_2, \xi_1 +\xi_2, z_1 +z_2 +\frac{\xi_1 s_2 -\xi_2 s_1}{2})
\end{equation}

The Schroedinger representation is the unitary representation $\varrho$ of $H$ on $L^2(\mathbb{R})$ defined as
\begin{equation}
\varrho_x f(t) = e^{iz + i\xi (t-\frac{s}{2})} f(t-s)
\end{equation}
The short time Fourier transform with window h can then be described as
\begin{equation}
Sf(t,\omega) = (f, \varrho_{(t,\omega, 0)}h)
\end{equation}
with the hermitian scalar product $(.,.)$ in $L^2(\mathbb{R})$.

The window $h$ usually is taken as a smooth function localized near the origin, that decreases quickly for big values of the argument. A typical example is the Gaussian distribution
$g(t) = \frac{1}{\sqrt{2 \pi}}e^{-\frac{t^2}{2}}$. In this case, the transfomation is often called the Gabor transformation, with reference to the paper by Gabor from 1946.
The function $f \in L^2(\mathbb{R})$ can be recovered from $Sf$ as

\begin{eqnarray}
2\pi \|h\|^2 f(t)&=&  \int\int Sf(s,\xi)  \varrho_{(s,\xi, 0)}h(t) ds \, d\xi
\end{eqnarray}
The Schroedinger representation induces an action on the short time Fourier transform. Setting $x=(t,\omega, 0)$ and $y=(s,\xi,0)$ formula (3) gives
\begin{eqnarray}
S\varrho_y f(t,\omega) &=& (\varrho_y f, \varrho_x h)=(f,\varrho_{y^{-1}.x}h)\\
&=& e^{i\frac{\xi t-\omega s}{2}}Sf(t-s,\omega -\xi)
\end{eqnarray}
since
\begin{equation}
y^{-1}.x= (t-s,\omega -\xi, \frac{s\omega -t\xi}{2})
\end{equation}

The action
\begin{eqnarray}
S\varrho_{(s,\xi,z)}^2 F(t,\omega)
:&=& e^{iz + i\frac{\xi t-\omega s}{2}}F(t-s,\omega -\xi)
\end{eqnarray}
is unitary on $L^2(\mathbb{R}^2)$ and the short time Fourier transform intertwines the Schroedinger representation with it:
\begin{equation}
S\varrho_x f = \varrho_x^2 (S f)
\end{equation}
for $x \in H$.
When reconstructing $f$ from its short time Fourier transform $Sf$, the contribution of $\varrho_{(t,\omega,0)}h$ is weighed with the coefficient $Sf(t,\omega)$. If in the rearrangement process the coefficient is considered at the point $(\tilde{t},\tilde{\omega})$ instead of at the point $(t,\omega)$, then a phase change has to be expected. It is determined by the group action and the window $h$:

\begin{eqnarray}
arg(\varrho_{(\tilde{t},\tilde{\omega},0)} h ,\varrho_{(t,\omega,0)}h)&-&arg(\varrho_{(t,\omega,0)}h,\varrho_{(t,\omega,0)}h)\\&=&  arg(\varrho_{(\tilde{t},\tilde{\omega}, 0)}h ,\varrho_{(t,\omega, 0)}h)\\
 &=& arg(h,\varrho_{(-\tilde{t},\tilde{-\omega},0).(t,\omega, 0)}h)\\
&=&\frac{\tilde{\omega} t-\tilde{t}\omega}{2} + arg\, Sh(t-\tilde{t},\omega -\tilde{\omega})
\end{eqnarray}

As an example, take the window to be the normalized Gaussian function
\begin{equation}
g(t)=\frac{1}{\sqrt{2 \pi}} e^{-\frac{t^2}{2}}
\end{equation}
Then
\begin{equation}
Sg(t,\omega)= e^{-\frac{t^2 + \omega^2}{2}}
\end{equation}
such that the expected change in the argument is just $\frac{\tilde{\omega} t-\tilde{t}\omega}{2}$.
Observe that the Gaussian function is an extremal for the Heisenberg uncertainty inequality (i.e. for this function there is equality in the uncertainty relation). The function is also special in that $Sf(t,\omega)$, calculated with the window $g$ is up to the factor $e^{-\frac{t^2 +\omega^2}{2}}$ a holomorhic function in $z=\omega +i t$:

\begin{eqnarray}
Sf(t,\omega)&=&\frac{1}{\sqrt{2 \pi}} \int e^{-\frac{(s-t)^2}{2}} e^{i\omega (s-\frac{t}{2})} f(s) ds\\&=& e^{-|z|^2} e^{-\frac{(\omega +it)^2}{4}} \frac{1}{\sqrt{2 \pi}}\int f(s) e^{-\frac{(s-t+i\omega)^2}{2}}
\end{eqnarray}

\section{Reassignment for the spectogram}

Reassignment is a method of postprocessing spectrograms with the aim to bring forward the structural components that may be hidden in it. The essential component of reassignment is a parameter transformation $T_f: (t,\omega)\rightarrow (\tilde{t},\tilde{\omega})$ that depends on the spectogram itself, more precisely on its phase $\varphi$. The transformation $T=T_f$ is

\begin{eqnarray}
\tilde{t}&=&\frac{t}{2} -\frac{\partial \varphi}{\partial \omega}\\
\tilde{\omega}&=&\frac{\omega}{2} + \frac{\partial \varphi}{\partial t}
\end{eqnarray}
If it is injective, with $J_T (t,\omega)$ its Jacobian determinant different from $0$, then the transformed function is defined as
\begin{equation}
\tilde{Sf}(\tilde{t},\tilde{\omega}) = e^{-i\frac{\tilde{t}\omega -\tilde{\omega}t}{2}}Sh(\tilde{t}-t,\tilde{\omega}-\omega) Sf(t,\omega)|J_T(t,\omega)|^{-1}
\end{equation}
The phase factor is the expected phase change in the passage from $(t,\omega)$ to $(\tilde{t},\tilde{\omega})$ as calculated in the first section. If $T_f$ fails to be injective or if $J_T$ has zeroes, then $\tilde{Sf}$ can be defined as the (complex) measure such that for all Borel sets $A$
\begin{equation}
\int_A d\tilde{Sf}(\tilde{t},\tilde{\omega}) = \int_{T^{-1}A}e^{-i\frac{\tilde{t}\omega -\tilde{\omega}t}{2}}Sh(\tilde{t}-t,\tilde{\omega}-\omega) Sf(t,\omega)|J_T(t,\omega)|^{-1} dt\,d\omega
\end{equation}
The notation is consistent in that
\begin{equation}
d\tilde{Sf}(\tilde{t},\tilde{\omega}) = \tilde{Sf}(\tilde{t},\tilde{\omega}) d\tilde{t}\,d\tilde{\omega}
\end{equation}
for injective transformations.\\
In practice, $Sf$ and its reassignment are calculated on fine grids. $T_f$ transforms a grid $\{(t_j,\omega_k)\}$ of points in the $(t,\omega) $ plane into a set of points in the $(\tilde{t},\tilde{\omega})$ plane. The image points $T(t_j,\omega_k)$ that have $(\tilde{t}_l,\tilde{\omega}_m)$ as their closeset neighbour out of a grid of points in the $(\tilde{t},\tilde{\omega})$ plane contribute to the sum in the definition of the reassigned value of $Sf$ at the grid point $(\tilde{t}_l,\tilde{\omega}_m)$

\begin{equation}
\tilde{Sf}(\tilde{t}_l,\tilde{\omega}_m) =\sum e^{i\frac{\tilde{\omega_m}t_j - \tilde{t}_l \omega_k}{2}}Sf(t_j,\omega_k)
\end{equation}

In [CDAF] it is shown that the displacement vector
\begin{eqnarray}
v(t,\omega) &=& (\tilde{t}-t,\tilde{\omega}-\omega)\\
&=&(-\frac{t}{2} -\frac{\partial \varphi}{\partial \omega}(t,\omega),
-\frac{\omega}{2} + \frac{\partial \varphi}{\partial t}(t,\omega))
\end{eqnarray}
at any point $(t_0,\omega_0)$ is tangent to the level lines of the "geometric phase function" $\Phi_{(t_0,\omega_0)}$ associated with the point $(t_0,\omega_0)$:
\begin{eqnarray}
\Phi_{(t_0,\omega_0)}(t,\omega):&=& arg (\varrho_{(-t_0,\omega_0,0)}f,\varrho_{(t,\omega,0)}h)\\
&=& \varphi(t+t_0,\omega +\omega_0) -\frac{t\omega_0 -t_0\omega}{2}
\end{eqnarray}
Clearly at $(t_0,\omega_0)$
\begin{eqnarray}
v(t_0,\omega_0) &=&(-\frac{t_0}{2} -\frac{\partial \varphi}{\partial \omega},
-\frac{\omega_0}{2} + \frac{\partial \varphi}{\partial t})\\
&=& (-\frac{\partial \Phi_{(t_0,\omega_0)}}{\partial \omega}, \frac{\partial \Phi_{(t_0,\omega _0})}{\partial t})|_{(t,\omega)=(0,0)}
\end{eqnarray}
is orthogonal to $grad\,\Phi_{(t_0,\omega_0)}$.\\
It furthermore was observed in [CDAF] that the displacement vector can be derived
from the function $F(z,\bar{z})$ in the variable $z=\omega +it$ and its conjugate that is defined as
\begin{equation}
Sf(t,\omega) = F(z,\bar{z})\,e^{-\frac{|z|^2}{4}}
\end{equation}
In particular, if the window $h$ is taken to be the Gaussian function, then $F(z)$ is holomorphic (i.e. independent of $\bar{z}$ ) and
\begin{equation}
v(t,\omega) = grad\,log\,|F(z)|
\end{equation}

\section{The wavelet transform}

The affine group
\begin{eqnarray*}
G &=& \{(a,b): \, a>0, b \in \mathbb{R}\}\\
&=& \{(a,b): \, t \rightarrow \, at+b,\, a>0, \, b \in \mathbb{R}\}
\end{eqnarray*}

with group law
\begin{eqnarray*}
(a,b)(a',b')=(aa',ab'+b)
\end{eqnarray*}
acts on $L^2 (\mathbb{R})$. Translations and dilations induce actions
on $L^2(\mathbb{R})$. The translations
$$\tau_b: t \rightarrow t+b$$
describe the time shift. This action is quite natural and reflects
time invariance. The dilations

$$\delta_a: \; t \, \rightarrow\, at$$
however enter the picture only at the level of the hearing
process. They reflect that the ear recognizes a sound
independently of its basic pitch.\\

The induced actions on $L^2(\mathbb{R})$ are denoted with the same
symbols:
\begin{eqnarray*}
  \delta _a f(t) &=& \frac{1}{\sqrt{a}}f(\frac{t}{a}) \\
  \tau_b f(t) &=& f(t-b)
\end{eqnarray*}
The factor $\frac{1}{\sqrt{a}}$ is introduced such that the
action is unitary. We then set
$$\varrho_{a,b}f(t)=\frac{1}{\sqrt{a}} f(\frac{t-b}{a}) = \tau_b
\delta_a f(t)$$

Under this action, the subspace $L^2(\mathbb{R},\mathbb{R})$ of real valued
functions is preserved.
Under the Fourier transform
\begin{equation}
    f(\omega) = \frac{1}{\sqrt{2 \pi}}\int e^{-it\omega} f(t) \,dt
\end{equation}
the action is intertwined to the action

\begin{equation}
    \hat{\varrho}(a,b)f(t)=e^{-it\omega}\sqrt{a}\hat{f}(a\omega)
\end{equation}

A wavelet is a function $\psi\in L^2(\mathbb{R})$ with mean value zeroe
and normalized by the condition $\|\psi\|=1$. We will in the following consider real valued wavelets and analytic wavelets. They have to satisfy the admissibility condition

\begin{equation}
    \int_0^{\infty}\frac{|\hat{\psi}(\omega)|}{\omega} \,d\omega = C <\infty
\end{equation}
The wavelet $\psi$ is analytic, if $\hat{\psi}(\omega)=0$ for all $\omega <0$

 Under the action of the affine group $G$ the wavelet is transformed into the functions
\begin{equation}
    \psi_{a,b}=\varrho (a,b)\psi
\end{equation}
with shifted mean localization and frequency content. Since $\varrho$ is
unitary, the functions $\psi_{a,b}$ are still normalized. The
wavelet transform of $f\in L^2(\mathbb{R})$ at scale $a$ and
time $b$ is
\begin{equation}
    W_\psi f(a,b)=(f,\psi_{a,b})=\int_{-\infty}^{\infty}f(t)
    \frac{1}{\sqrt{a}}\bar{\psi}(\frac{t-b}{a})\,dt
\end{equation}
The function $f$ can be reconstructed from its wavelet transform

\begin{equation}
f(t) = \frac{1}{C} \int_0^{\infty}\int_{-\infty}^{\infty}\frac{da \, db}{a^2}    W_\psi f(a,b) \varrho_{a,b}\psi(t)
\end{equation}

The density $Wf(a,t)\frac{da\, db}{a^2}$ can be considered as the weight with which the function $\varrho_{a,b}\psi$ appears.

The wavelet transform intertwines the action $\varrho$ with the action $\varrho^2$ on $L^2(\mathbb{R}^2)$:

\begin{eqnarray*}
  W \varrho_{c,s}f(t) &=& (\varrho_{c,s}f,\varrho_{a,b}\psi) = (f,\varrho_{(c,s)^{-1}(a,b)}\psi)\\
  &=& W f((c,s){-1}(a,b))\\
  :&=& \varrho_{c,s}^2 Wf(a,b)
\end{eqnarray*}
The representation $\varrho^2$ of $G$ is unitary on $L^2( (0,  \infty)\times \mathbb{R},\frac{da\, db}{a^2})$. The measure $\frac{da\, db}{a^2}$ is the Haar measure for the group action.

\section{Reassignment for the scalogram}

The scalogram for the wavelet transform of the function f is the graph of $Wf(a,b)$. It will be assumed that $f$ is real valued. If the wavelet is also real valued, then the scalogram is real valued. However if the wavelet is analytic, then $Wf$ is complex valued and its phase is well defined as long as $|Wf(a,b)| \neq 0$. In this case, the transformation $T$ in the rearrangement process is defined as

\begin{eqnarray*}
\frac{1}{\tilde{a}}&=& \frac{\partial \varphi}{\partial t}(a,t)\\
\tilde{t}&=&t + a^2 \frac{\partial \varphi}{\partial a}(a,t)
\end{eqnarray*}
  Phase and amplitude changes have to be expected when changing $(a,t)$ to $ (\tilde{a},\tilde{t})$. They are given by

\begin{equation}
( \varrho_{\tilde{a},\tilde{t}}\psi,\varrho_{a,t}\psi)= W\psi((\tilde{a},\tilde{t})^{-1}(a,t)) =W\psi(\frac{a}{\tilde{a}},\frac{b-\tilde{b}}{\tilde{a}})
\end{equation}

The transformation of the spectogram is then defined by

\begin{equation}
\tilde{W}f(\tilde{a},\tilde{t}) = ( \varrho_{\tilde{a},\tilde{t}}\psi,\varrho_{a,t}\psi) Wf(a,t)\, |J_T(a,t)|^{-1}
\end{equation}
The results in section 3 show that
the reassignment method for the spectogram works particularly well if the window is a Gaussian, i.e. if it is an extremal for the Heisenberg uncertainty relation. Similarly, reassignment of the scalogram is particularly suited for wavelets derived from extremals of the uncertainty relation for the affine group. The Fourier transforms of the extremals as given in
[R] are

\begin{equation}\hat{h}(\omega) = k e^{i\varepsilon sgn(\omega) -i\alpha sgn(\omega)\,log |\omega| -i\beta
\omega}\,e^{-\kappa |\omega|}\,|\omega|^{\kappa \nu
-\frac{1}{2}}\end{equation} with real constants $k, \varepsilon, \alpha,
\beta, \kappa $ and $\nu $. Square integrability implies $\kappa > 0$.
From the explicit form it is clear that the space of extremals is
invariant under the action of $G$.\\
In the following, analytic wavelets of the form $\psi = h +i Hh$ with h an extremal and H the Hilbert transform are considered. The Fourier transform of $Hh$ is $-isgn(\omega)\hat{h}(\omega)$
and hence
$$\hat{\psi}(\omega) =
\begin{cases}
2\hat{h}(\omega) & \omega > 0\\
0  & otherwise
\end{cases}
$$
In this situation, the wavelet transform with wavelet $\psi$ satisfies the differential equation

\begin{equation}
 a \frac{\partial }{\partial a}log\,Wf =
\kappa  -i\alpha -a(\beta
-i\kappa)\frac{\partial}{\partial t}log\,Wf
\end{equation}\\
(see [R], equation 62). The wavelet transform is written in polar coordinates
$$Wf(a,t) = r(a,t)e^{i\varphi(a,t)}$$
In the case $\beta =0$ the real and imaginary parts of the equation give

\begin{eqnarray}
a\frac{\partial\varphi}{\partial a} &=& -\alpha +a\kappa  \frac{\partial}{\partial t} log\,r\\
a\frac{\partial \varphi}{\partial t} &=& 1 - \frac{a}{\kappa}\frac{\partial}{\partial a}log\,r
\end{eqnarray}

This shows, that the rearrangement transformation T can also be defined in terms of the logarithmic amplitude:
\begin{eqnarray*}
\frac{1}{\tilde{a}}&=&  \frac{\partial \varphi}{\partial t}(a,t) = \frac{1}{a} - \frac{1}{\kappa}\frac{\partial}{\partial a}log\,r\\
\tilde{t}&=&t + a^2 \frac{\partial \varphi}{\partial a}(a,t) = t -a \alpha +a^2 \kappa  \frac{\partial}{\partial t} log\,r
\end{eqnarray*}

In the case $\beta \neq 0$, the reassignment should be defined by the transformation $T_{\beta}$:

\begin{eqnarray*}
\frac{1}{\tilde{a}}&=& \frac{\partial \varphi}{\partial t}(a,t)\\
\tilde{t}&=&t + a^2 (\frac{\partial}{\partial a} + \beta
\frac{\partial}{\partial t})\varphi (a,t)
\end{eqnarray*}
This choice is suggested when one considers the new time parameter $s=t-a \beta$ that takes a shift proportional to scale into account (intuitively, there is a natural time delay when the pressure wave travels along the cochlea). If one would set
\begin{equation}
r(a,t) e^{i\varphi (a,t)} =R(a,s)e^{i\Phi(a,s)}
\end{equation}
then the reassignment transformation would have to be based on

\begin{eqnarray}
\frac{\partial\Phi}{\partial a} &=&   \frac{\partial \varphi}{\partial a} +\frac{\partial \varphi}{\partial t}\frac{\partial t}{\partial a} =\frac{\partial \varphi}{\partial a} + \beta \frac{\partial \varphi}{\partial t}\\
\frac{\partial \Phi}{\partial s} &=& \frac{\partial \varphi}{\partial t}\frac{\partial t}{\partial s} =\frac{\partial \varphi}{\partial t}
\end{eqnarray}
This then provides the motivation for the reassignment transformation $T_\beta$.\\

In [R] the uncertainty principles for the affine group are discussed for general parameters. The Fourier transforms of the extremals for the general uncertainty principle are given as

\begin{equation}
\hat{h}_c(\omega) = k e^{i\varepsilon sgn(\omega) -i\alpha
\,sgn(\omega)\,log |\omega| -i\beta \omega}\,e^{-\frac{\kappa}{c}
|\omega|^c}\,|\omega|^{ \kappa \nu -\frac{1}{2}}
\end{equation}
 with real
constants $k, \varepsilon, \alpha,\beta, \kappa$ and $\nu$.
 These solutions are in $L^2$ if both $\kappa > 0$ and $\nu > 0$.\\

It is further shown that the wavelet transforms with analytic wavelet $\psi = h_c + i H h_c$ approximately satisfy the differential equation

\begin{equation}
 a \frac{\partial }{\partial a}log\,Wf \cong \gamma
 -i\alpha -a(\beta -i\gamma)\frac{\partial}{\partial
t}log\,Wf
\end{equation}\\
Here, the constants are $\varepsilon =0$, $\nu =1$ and $\gamma = c \kappa$. (Observe that the function
$Zf$ in [R] is equal to the present function $Wf$ divided by $\sqrt{a}$).
This equation is referred to as the structure equation (up to this little modification).

The differential equation

\begin{equation} a \frac{\partial }{\partial a}Y(a,t) =
\gamma  -i\alpha -a(\beta
-i\gamma)\frac{\partial}{\partial t}Y(a,t)\end{equation}

has the particular solution

\begin{equation}P_{\gamma}(a) = (\gamma  - i \alpha )log\,a
\end{equation}
Its distinguished feature is the time independence. The
homogeneous equation associated to the (complex) structure
equation is
\begin{equation}
\frac{\partial }{\partial a} Y(a,t) = -(\beta
-i\gamma)\frac{\partial}{\partial t}Y(a,t)
\end{equation}

With the variable change
\begin{eqnarray*}
z &=& t-a\beta +ia\gamma \\
\bar{z} &=&  t-a\beta -ia\gamma
\end{eqnarray*}

this becomes a $\bar{\partial}$-equation in z. The solutions are holomorphic functions in the variabble $z$. The solutions of the linear inhomogeneous
equation have the representation
\begin{equation}Y(a,t) = P_{\gamma}(a) + F(z) \end{equation}
with $F$ a holomorphic function in the variable $z = t-a\beta +ia\gamma $. Since $a > 0$
(and $\kappa > 0$) it is defined in the upper half space $\{ z\in
\mathbb{C}:\,Im\,
z>0\}$. \\

Real and imaginary part of the structure equation give

\begin{eqnarray*}
a(\frac{\partial}{\partial a}+\beta\frac{\partial}{\partial t}) log\,r &\cong& \gamma - a\gamma \frac{\partial \varphi}{\partial t}\\
a(\frac{\partial}{\partial a}+\beta\frac{\partial}{\partial t})\varphi  &\cong& -\alpha - a\gamma \frac{\partial}{\partial t} log\,r
\end{eqnarray*}
The reassignment
$T_\beta$ is therefore

\begin{eqnarray*}
\frac{1}{\tilde{a}}&=& \frac{\partial \varphi}{\partial t}(a,t) \cong \frac{1}{a}-\frac{1}{\gamma} (\frac{\partial}{\partial a} + \beta
\frac{\partial}{\partial t}) log\,r(a,t)\\
\tilde{t}&=&t + a^2 (\frac{\partial}{\partial a} + \beta
\frac{\partial}{\partial t})\varphi (a,t) \cong t -a \alpha -a^2 \gamma
\frac{\partial}{\partial t} log\, r (a,t)
\end{eqnarray*}

The reassignment can be expressed in terms of the holomorphic function $F$.
Since $z = t - a\beta + ia\gamma$ one obtains

\begin{equation}
\frac{\partial F}{\partial a} + \beta \frac{\partial F}{\partial t} = i \gamma \frac{\partial F}{\partial z}
\end{equation}
Starting from the approximate equation
\begin{equation}
log\,Wf(a,t) \cong P_{\gamma} (a) + F(z)
\end{equation}
this leads to

\begin{eqnarray*}
a(\frac{\partial}{\partial a} + \beta
\frac{\partial}{\partial t}) log\,r (a,t) &\cong& \gamma -\gamma a Im
\frac{\partial F}{\partial z}\\
a(\frac{\partial}{\partial a} + \beta
\frac{\partial}{\partial t})\varphi (a,t) &\cong& - \alpha +\gamma a Re
\frac{\partial F}{\partial z}
\end{eqnarray*}

\begin{eqnarray*}
\frac{1}{\tilde{a}} &\cong&
Im \frac{\partial F}{\partial z} \\
\tilde{t}&\cong& t -a \alpha + a^2 \gamma
Re \frac{\partial F}{\partial z}
\end{eqnarray*}

\section{Examples}
1)  The holomorphic function associated to the signal $f(t)= cos \,\nu\,t$ under the analytic wavelet transformation with the wavelet $\psi = h_{\gamma} + i Hh_{\gamma}$ is $F(z) = const +i\nu z$. The reassignment transformation $T_{\beta}$ therefore is

\begin{eqnarray*}
\frac{1}{\tilde{a}} &=& \nu\\
\tilde{t}&=& t -a \alpha
\end{eqnarray*}

2) The holomorphic function associated to a click (the $\delta$-function at $0$) under the same wavelet transformation as above is
$$ F(z)= const + (\gamma -i\alpha)log\,z$$
Since
\begin{equation}
\frac{\partial F}{\partial z} =\frac{\gamma -i\alpha}{z} =\frac{(\gamma -i\alpha)(t-a\beta +ia\gamma)}{|z|^2}
\end{equation}

the reassignment transformation is

\begin{eqnarray*}
\frac{1}{\tilde{a}} &\cong&
 \frac{a\gamma^2  - \alpha (t-a\beta) }{|z|^2 }\\
\tilde{t}&\cong& t -a \alpha + a^2 \gamma ^2
 \frac{(t-a\beta) + a\alpha }{|z|^2}
\end{eqnarray*}
with $|z|^2=(t-a\beta)^2 +a^2 \gamma ^2$.
At $t-a\beta =0$ this is

\begin{eqnarray*}
\frac{1}{\tilde{a}} &\cong&
 \frac{1}{a}\\
\tilde{t}&\cong& t
\end{eqnarray*}


\vspace{4cm}The author would like to acknowledge the support of the  Centro di Ricerca Matematica Ennio De Giorgi in Pisa

\vspace{4cm} H.M.Reimann, Institute of Mathematics, University of
Berne, Sidlerstrasse 5, 3012 Berne, Switzerland

\end{document}